\documentclass{article}
\usepackage{amssymb,amsthm,amsmath,amsfonts}
\usepackage{epsfig}

\usepackage[a4paper,bindingoffset=0.2in,%
            left=.5in,right=.5in,%
            footskip=.25in, top = .5in]{geometry}

\theoremstyle{plain}
\begingroup
 
\newtheorem{theorem}{Theorem}

\endgroup
\newcommand{\nwc}{\newcommand}
\nwc{\Levy}{L\'{e}vy}
\nwc{\Holder}{H\"{o}lder}
\nwc{\cadlag}{c\`{a}dl\`{a}g}
\nwc{\be}{\begin{equation}}
\nwc{\ee}{\end{equation}}
\nwc{\ba}{\begin{eqnarray}}
\nwc{\ea}{\end{eqnarray}}
\nwc{\la}{\label}
\nwc{\nn}{\nonumber}
\nwc{\Z}{\mathbb{Z}}
\nwc{\C}{\mathbb{C}}
\nwc{\E}{\mathbb{E}}
\nwc{\R}{\mathbb{R}}
\nwc{\N}{\mathbb{N}}
\nwc{\prob}{\mathbb{P}}
\nwc{\Skor}{\mathbb{D}}
\nwc{\PP}{\mathcal{P}}
\nwc{\M}{\mathcal{M}}
\nwc{\law}{\stackrel{\mathcal{L}}{\rightarrow}}
\nwc{\eqd}{\stackrel{\mathcal{L}}{=}}
\nwc{\vp}{\varphi}
\nwc{\Vp}{\Phi}
\nwc{\veps}{\varepsilon}
\nwc{\eps}{\ve}
\nwc{\qref}[1]{(\ref{#1})}
\nwc{\D}{\partial}
\nwc{\dnto}{\downarrow}
\nwc{\nsup}{^{(n)}}
\nwc{\ksup}{^{(k)}}
\nwc{\jsup}{^{(j)}}
\nwc{\nksup}{^{(n_k)}}
\nwc{\inv}{^{-1}}
\nwc{\one}{\mathbf{1}}
\nwc{\argmin}{\mathrm{arg}^+\mathrm{min}}
\nwc{\argmax}{\mathrm{arg}^+\mathrm{max}}
\nwc{\Rplus}{\R_+}
\nwc{\Rorder}{\R_<}
\nwc{\xx}{\mathbf{x}}
\nwc{\emp}{\mu}
\nwc{\empN}{F^n}
\nwc{\lossN}{L^n}
\nwc{\Lip}{\mathrm{Lip}}
\nwc{\BL}{\mathrm{BL}}
\nwc{\ddm}{d}
\nwc{\no}{\textcolor{red}{[CLUNK!]}}

\theoremstyle{definition}
\newtheorem{defn}[theorem]{Definition} 
\newtheorem{remark}[theorem]{Remark}

\theoremstyle{remark}
 
\numberwithin{equation}{section}
\numberwithin{figure}{section}

\usepackage{color}

\begin{document}

\title{Renewal-scaled solutions of the Kolmogorov forward equation for residual times \\}

\author{Joseph Klobusicky\textsuperscript{1}}

\date{}

\maketitle

\begin{abstract}
Let $N(\tau)$ be a renewal process for independent holding times $\{X_i\}_{k \ge 0}$  ,where $\{X_k\}_{k\ge 1}$ are identically distributed with   density $p(x)$.  If the associated residual time $R(\tau)$ has a density $u(x,\tau)$, its Kolmogorov forward equation  is given by \be \nn
\partial_\tau u(x,\tau)-\partial_x u(x,\tau) = p(x)u(0,\tau), \quad x,\tau \in [0, \infty),
\ee
with an initial holding time density $u(x,0)=u_0(x)$.  We derive a measure-valued  solution formula for the density of residuals times  after an expected number of renewals occur. Solutions under this time scale are then shown  to evolve continuously in the space of measures with the weak topology under mild conditions.  \end{abstract}



\medskip
\noindent
\footnotetext[1]
{Department of Mathematics, The University of Scranton, 800 Linden St., Scranton, PA 18510.
Email: joseph.klobusicky@scranton.edu.}
 \section{Introduction}

Let $\{X_i\}_{i\ge 1}$ be a sequence of iid  random variables in $\mathbb R_+ = [0,\infty)$ called \textit{holding times}, each having a  density $p(x)\in L^1(\mathbb R_+)$.  Also define a random initial  holding time $X_0 \in \mathbb
R_+ $,  independent from $\{X_i\}_{i\ge 1}$, but possibly distinct in law. The renewal process $N(\tau)$ associated with $\{X_i\}_{i\ge 0}$ counts the number of renewals up to time $\tau$, and is given by  
\begin{equation}
N(\tau) = \sup\left\{i\ge 0:\sum_{j = 0}^i X_i\le \tau \right\}. \label{renewal}
\end{equation}
Related to the renewal process is the \textit{residual time} (also known as the forward recurrence time), defined as     
\begin{equation}
R(\tau) = \sum_{i = 1}^{N(\tau)+1} X_i-\tau. \label{residual} 
\end{equation}
At time $\tau$,  $R(\tau)$  is the time remaining until the next renewal.  Residual times are  fundamental  in  renewal theory and queueing theory (see  Ch. 9 of Cox \cite{cox2017theory}  for a detailed introduction), and may be viewed as an extension of homogeneous Poisson processes, in which holding times, and consequently residual times, are exponentially distributed \cite{resnick2013adventures}. 

If  $R(\tau)$ has a differentiable density in time and space   $u(x,\tau)\in C^1(\mathbb R^+ \times \mathbb R^+, \mathbb
R^+)$, then  the  Kolmogorov  forward equation for residual times is given by
 \begin{align}  \label{pde}
\partial_\tau u(x,\tau) - \partial_x u(x,\tau) &= p(x)u(0,\tau), \quad \tau,
x \ge0, \quad u(x,0) = u_0(x)
.  
\end{align}
     This equation has  the explicit solution (see Pg. 63 of \cite{cox1967renewal}) of 
\begin{equation}\label{intsol}
u(x,\tau) = u_0(x+\tau) +\int_0^\tau p(x+\tau-s)\alpha(s)ds. 
\end{equation}
Here,  $\alpha(s)$ is the renewal density
\begin{equation}
  \alpha(\tau)  = u_0(\tau) *\sum_{i = 0}^\infty p^{*(i)}(\tau),\label{trace}
\end{equation}
where $p^{*(i)}$ denotes $i$-fold self-convolution (with the convention $p^{*(0)}
= 1$ and $p^{*(1)} = p$).
An important fact  (see Sect. 9.12 in \cite{stewart2009probability})  is that
\begin{equation}
A(\tau):= \mathbb E[N(\tau)]= \int_0^\tau \alpha(s)ds.
\end{equation}
In Section \ref{subsec:strong} we derive  (\ref{pde})
 through asymptotic expansions, and in Section \ref{subsec:classic},  we derive (1.4)  through Laplace transforms,
 a method similar to     \cite{kim2004markovian}.

The main novelty of this paper is a change of variables for (\ref{pde}) based on expected renewals, presented in Section \ref{meassect}. Specifically, we introduce the new time scale $t(\tau) = A(\tau).$ Using the $t$ time scale is  a natural option  for
several scenarios in queuing theory. For an example, consider a factory lit by a large collection
 of  $N$ light bulbs,  running simultaneously.  When a bulb breaks, it is immediately
replaced by a new light bulb  with a random run length distributed with
respect to $p(x) $. If the initial distribution of run lengths at time zero for the $N$ bulbs   is approximately
$u_0(x)$, then   $u(x,\tau)$ estimates remaining run lengths    at time  $\tau$. What,
then, is the distribution of remaining run lengths  after    $N t$  light bulbs
 have been replaced? The answer to this question is   provided through the  renewal-scaled solution formula (\ref{measform}). Given the   expansive history of renewal theory, the author   expected  to  find (\ref{measform})  in the literature, but was unable to do so after a thorough search.      

Both the change of time scale and method for renewal-scaled solutions  are motivated by  the work of  Menon,
Niethammer, and Pego in \cite{menon2010dynamics}, who investigated a wide range of clustering events generalizing
the 1D Allen-Cahn equation in mathematical physics.
 The measure-valued space of clusters was shown to be continuous
in time in the space of probability measures $\mathcal P(\mathbb R^+)$ through  using an intrinsic
time scale based on the number of clusters in the system.  In our case, we  define measure-valued solutions of (\ref{pde}) in the $t$ time scale which evolve continuously in  $\mathcal P(\mathbb R^+)$  under a wide range of holding times. We conclude with an example in  Section \ref{ptmass}, which illustrates  with  how renewal-scaled solutions continuously evolve under initial holding times with  point masses.

\section{Strong solutions} \label{sec:strong}

 \subsection{Derivation of forward equation} \label{subsec:strong}

\begin{figure}
\begin{centering}
\includegraphics[width=.7\textwidth]{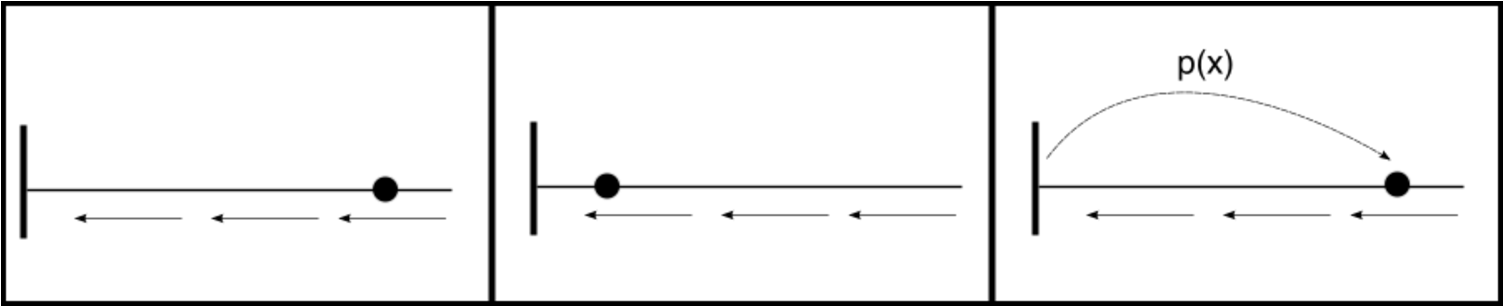}
\caption{\textbf{A dynamical system representation} Left: A particle  at time $\tau =
0$ is placed on the positive real line. Middle: The particle moves toward
the origin  at unit speed. Right: The particle reaches the origin, where
where it is randomly reassigned to $(0,\infty)$ under
a probability density $p(x)$.} \label{figure}
\end{centering}
\end{figure}

In this section will use  an equivalent reformulation for $R(\tau)$  from  the perspective
of a random dynamical system.
This is done by considering a particle at time
$\tau = 0$ randomly placed on  $\mathbb R^+ $,
with its initial position $R(0)$ given by the first holding time $X_0$. The
particle
drifts toward the origin at unit speed until it reaches the origin at time
$\tau = X_0$. At this time, it  is immediately reassigned to $R(X_0) = X_1
\in \mathbb R^+$, distributed with respect to  $p(x)$.
The particle proceeds as before, moving toward the origin and being reassigned
upon its arrival, with  $R(\sum_{j = 1}^{k-1} X_{j}) = X_{k}, $ and   $\{X_i\}_{i\ge
1}$ iid (see Figure \ref{figure}).
 It is clear that the position of the particle at time $\tau$ corresponds
with the residual time $R(\tau)$ under holding times $\{X_k\}_{k\ge 1}.$

To derive (\ref{pde}), we will compute asymptotics for probabilities that a particle is in a small interval.  Assume for now that $p(x)$ is continuous, and $u(x,\tau)$
is differentiable in both $x$ and $\tau$. 
Then an  asymptotic expansion for the  probability that a particle is in $[x,x+\Delta x]$
at time $\tau$ is
\begin{equation} \label{pr}
\mathbb P(R(\tau) \in [x,x+\Delta x]) = u(x,\tau)\Delta x +\frac 12\partial_x
u(x,\tau)(\Delta
x)^2+o((\Delta x)^2). 
\end{equation}
Next, we  consider the probability that the particle  at initial time $\tau-\Delta
x$  is contained in $[x,x+\Delta x]$ at time $\tau$.  This can occur in two
ways.  The simplest case is when the particle is initially in $[x+\Delta
x, x+2\Delta x]$. Denote  this event $A$, whose probability has the expansion
\begin{align}
\mathbb P(A) &= \mathbb P(R(\tau-{\Delta x}) \in [x+\Delta
x, x+2\Delta x])= u(x,\tau)\Delta x -\partial_\tau
u(x,\tau)(\Delta
x)^2+\frac 32\partial_x u(x,t)(\Delta
x)^2 +o((\Delta x)^2). \label{pa}
\end{align}

  The second case is when a particle begins  in $[0,\Delta x]$, reaches the origin, and is reassigned
 before drifting into $[x,x+\Delta x]$ at time
$\tau$. Denote this event $B$. We may partition this event according to the particle's initial
position $s \in [0,\Delta x]$, and how many times the particle jumps.
 Then
\begin{align} 
\mathbb P(B) &= \int_0^{\Delta x}u(s,\tau-\Delta x)[( \mathbb P((\hbox{1
jump}) \cap B| R(\tau-\Delta x)= s)  + \mathbb P((>\hbox{1
jump}) \cap
B| R(\tau-\Delta x) = s)]ds \nonumber\\
 &= \int_0^{\Delta x}u(s,\tau-\Delta x)\left(\int_{x+ \Delta x-s}^{x+2\Delta
x-s}p(y)dy\right)ds+o\big((\Delta x)^2 \big)   = u(0,\tau)p(x)(\Delta x)^2 +o\big((\Delta x)^2 \big). \label{pb}
\end{align}   
Since $\mathbb P(R(\tau) \in [x,x+\Delta x]) = \mathbb{P}(A)+ \mathbb{P}(B)$,
we  compare $(\Delta x)^2$ terms in (\ref{pr}), (\ref{pa}), and (\ref{pb}) to arrive at   (\ref{pde}).
\subsection{Classical solutions}\label{subsec:classic}

In this section, we will derive the solution formula (\ref{intsol}) for initial conditions in the space
\be
\mathcal C = \{u \in C^1(\mathbb R^+) \cap L^1(\mathbb R^+): u>0, \|u\|_1 = 1, u(x)\rightarrow 0 \hbox{ as } x\rightarrow \infty\}.
\ee
\begin{defn}
Let  $p(x) \in \mathcal C$. A function $u(x,\tau)  \in C^1(\mathbb R^+\times \mathbb
R^+; \mathbb R^+)$ with $u(x,0) \in \mathcal C$ which solves $(\ref{pde})$ is a \textbf{strong solution}. 
\end{defn}Our first step in finding an explicit strong solution formula is integrating (\ref{pde}) along characteristics $x(\tau) = x_0-\tau$,
which gives us the integral form of (\ref{pde})
\begin{equation}\label{intsolz}
u(x,\tau) = u_0(x+\tau) +\int_0^\tau p(x+\tau-s)u(0,s)ds. 
\end{equation}
By comparing (\ref{intsol}) and (\ref{intsolz}), we see that finding an explicit solution formula is equivalent to showing $u(0,\tau) = \alpha(\tau) $. This in fact follows immediately by simply setting $x = 0$ in (\ref{intsolz}), from which we obtain the well-known renewal equation whose solution is $\alpha(\tau)$. However, we derive here the explicit formula for $\alpha(\tau)$, as some of the intermediate formulas produced will be of use in later sections.  

Toward this end, we apply  the Laplace transform  
\begin{equation} \bar u(q,\tau) := \int_{\mathbb{R^+}}e^{-qx}u(x,\tau)dx, \quad q \in \mathbb C_+, 
\end{equation}
 to both sides of (\ref{pde}) to yield
\begin{equation} \label{lapform}
\frac{d}{d\tau} \bar u -q\bar u = (\bar p(q)-1)u(0,s).
\end{equation}
This may easily be solved for $\bar u$, with
\begin{equation}\label{lapsol}
\bar u(q,\tau) = e^{q\tau}\bar u_0(q)+(\bar p(q)-1)\int_0^\tau e^{q(\tau-s)}u(0,s)ds.
\end{equation}
This solution is similar to (\ref{intsol}), but now we can use Laplace inversion
to extract a formula for $u(0,s)$. Since we are taking the Laplace transform
of a probability density,  $|\bar u(q,\tau)|<1$ for $\tau>0, q \in
\mathbb{C}_+$. Thus $\bar u (q,\tau)e^{-q\tau} \rightarrow 0$ as $\tau \rightarrow
\infty$.  We then can obtain, from the $\tau\rightarrow \infty$ limit of
(\ref{lapsol}),
\begin{equation}\label{fpreinv}
\frac{1}{1-\bar p(q)}\bar u_0(q) = \int_{\mathbb{R^+}}e^{-qs}\alpha(s)ds  = [\bar u(0,\cdot)](q). 
\end{equation}
 We now have a formula for $\bar u(0,\tau)$ based on the initial
data.
Notice that the left hand side is a Laplace transform in the spatial variable,
whereas the right hand side is a Laplace transform  of $u(0,\tau)$ in the time variable.

To perform Laplace inversion on (\ref{fpreinv}), observe that $p$ is a probability density, so that  $|\bar p(q)|<1$ for $ q \in
\mathbb{C}_+$.
Thus, we can express the left hand side of (\ref{fpreinv}) as the geometric series 
\begin{equation}
 [\bar u(0,\cdot)](q) =\bar u_0(q) \sum_{i = 0}^\infty (\bar p(q))^i.  
\end{equation}
From the convolution formula, our renewal density has an explicit expression
\begin{equation}
  u(0,\tau)  = u_0(\tau) *\sum_{i = 0}^\infty p^{*(i)}(\tau),\label{trace}
\end{equation}
which is equivalent to $\alpha(\tau)$. Since $\alpha(t)$ is locally integrable (Th. 3.18 of \cite{liao2013applied}), (\ref{intsol}) is well-defined,  $u(x,\tau)$ is differentiable in both spatial and time variables, and for fixed $\tau \ge0$, $u(x,\tau)\rightarrow 0$ as $x\rightarrow \infty$.
Finally, we may show that $u(\tau,x)$ is a probability density for each $\tau\ge 0$  integrating (\ref{pde}) with respect to the $x$ variable. We summarize our findings in 

\begin{theorem}\label{strong} Let $p(x)\in C(\mathbb R^+),$ and $u_0(x)\in \mathcal C$.  Then there exists a unique strong solution $u(x,\tau)$ of (\ref{pde}) with $u(x,0) = u_0(x)$   given  by (\ref{intsol}). For each  $\tau>0$, $u(\cdot, t)\in \mathcal C$.  \end{theorem}

\section{Rescaling by expected renewals }\label{meassect}

In preparation for defining measure-valued solutions in the next section,  we'll first reformulate (\ref{pde}) in terms of  probability measures $\mathcal P(\mathbb R^+)$. In one dimension a probability  measure $\mu \in \mathcal P$ may may be identified with its cumulative distribution function $G_\mu(x) = \mu([0, x))$. In the future, when no confusion arises, we    will often write $G_\mu \in \mathcal P(\mathbb R^+).$ The Laplace transform of a measure $G$ is then defined as
\begin{equation}
\bar G(q)=   \int_0^\infty e^{-qx}
G(dx), \quad q \in \mathbb C^+. 
\end{equation}
If a measure $G(x)$ admits a density $G(dx) = g(x)dx$, then $\bar G(q) = \bar g(q)$. 

For $p(x),u_0(x) \in \mathcal C$,  (\ref{lapform}) may be rewritten in terms of measures as 
\begin{equation}\label{measode}
\frac d{d\tau}\bar U- qU = (1-\bar P(q))\alpha(\tau),
\end{equation}
where $dU_\tau(x) = u(x,\tau)dx$ and  $P(dx) = p(x)dx$.  Since $p$ and $u_0$ are strictly positive and continuous, so is $\alpha(\tau)$, and thus  $t = A(\tau)$ is strictly increasing and differentiable, with  $\frac{dt}{d\tau} = \alpha(\tau)$. We may now transform  (\ref{measode})  as 
\begin{equation}
 \frac 1{\alpha(\tau)}\frac d{d\tau}(e^{-q \tau}\bar U_\tau) =(1-\bar
P(q))e^{-q\tau}\nonumber.
\end{equation}
From the chain rule, we now convert to the $t$ time scale, with\begin{equation}
\frac d{dt}(e^{-q \tau(t)}\bar F_t) =(1-\bar P(q))e^{-q\tau(t)}, 
\end{equation}
where we now define   $ F_t(x) = U_{\tau(t)}(x)$. This in turn gives us the solution formula
\begin{equation}\label{solutionform}
\bar F_t(q)e^{-q\tau(t)} - \bar F_0(q) = (1-\bar P(q))\int_0^t e^{-q\tau(s)}ds.
\end{equation}
Since $\tau(t) \rightarrow \infty$ as $t\rightarrow \infty$, the limit of
(\ref{solutionform}) is then
\begin{equation}\label{initmeas}
\bar F_0(q)= (1-\bar P(q))\int_0^\infty e^{-q\tau(s)}ds,
\end{equation}
and consequently $\bar F_t$ takes the simple form
\begin{equation}\label{measform}
\bar F_t(q) = (1-\bar P(q))\int_t^\infty e^{q(\tau(t)-\tau(s))}ds.
\end{equation}

\subsection{Extension to measure valued solutions}

The solution formula (\ref{measform}) has only been defined for  the  strict class of functions $\mathcal C$, but its form suggests that we may  directly  extend solutions for measures   $F_0, P \in \mathcal P(\mathbb R^+)$. To do so,  however, we require a proper extension for $A(\tau)$, since the trace $\alpha(\tau) = u(0,\tau)$ assumes solutions have a density. Nonetheless,  we can use formula (\ref{trace}) to define 
\begin{equation}\label{measa}
A = F_0 *\sum_{i = 0}^\infty P^{*(i)},
\end{equation}
where convolution of the two measures $\mu$ and $\nu$ is defined as
\begin{equation}
\mu*\nu(E)  = \int_{\mathbb R^+} \int_{\mathbb R^+} \mathbf 1_{E}(x+y)d\mu(x)d\nu(y). \end{equation}
For a sequence of measures $f_n \rightarrow  \mu$, and $g_n \rightarrow  \nu$ converging weakly in $\mathbb R^+$, then by Prokhorov's theorem,  $f_n$ and $g_n$ are tight, and $f_n*g_n \rightarrow \mu *\nu$ weakly.
We may use this  to show that $A(\tau)$ is nondecreasing for arbitrary measures $F_0,P \in \mathcal P(\mathbb R^+)$ by regularization.

Even with a proper notion of $A(\tau)$,  still more is required to make 
(\ref{measform}) well-defined for measures.  This is because $A(\tau)$ may be constant   on an interval, or have jumps, both of which prohibit an inverse with a domain defined on all of  $\mathbb R^+$.  
To address this, we use a generalized  notion of inverse for nondecreasing functions.  For a nondecreasing function $ f:\mathbb R^+\rightarrow \mathbb R^+$, define the \textit{generalized inverse} $f^{\dagger}: \mathbb R^+\rightarrow \mathbb R^+$ by
\begin{equation}
  f^\dagger(q) =\inf \{s \in \mathbb R^+| f(s)>q\}. 
\end{equation}
When $f$ is strictly increasing and continuous, the usual inverse $f^{-1}$ and generalized inverse $f^{\dagger}$ coincide. For a generic distribution $A(\tau)$,  we  define
\be\label{meastau}
 \tau(t) = A^{\dagger}(t).
 \ee

 \begin{defn}
Let $P \in \mathcal P(\mathbb R^+)$.  A map $(F_0,t) \mapsto F_t \in \mathcal P(\mathbb R^+)$ defined through (\ref{measform}) is a \textbf{renewal-scaled solution} to (\ref{pde}), where $A(t)$ and $\tau(t)$ are defined by (\ref{measa}) and (\ref{meastau}).  
\end{defn}

\subsection{Properties of renewal-scaled solutions}

Here we will demonstrate that  the extension of the solution for  strong solutions to include measure-valued initial data and holding times is natural in the sense that measure-valued renewal-scaled solutions are weak limits of strong solutions in the $t$ time scale. We also show that our rescaling has the effect of smoothing  point masses approaching the origin, as weak solutions evolve continuously so long as $A^\dagger(t)$ is strictly increasing.
\begin{theorem}
Let $P \in \mathcal P(\mathbb{R}^+)$ be nondefective  (meaning $P(0) <1$). The following hold:

(1) Let $u^n(x,\tau)$ be a sequence of  strong solutions  with initial conditions $u_0^n(x)\in \mathcal C$ and holding time densities $p^n(x)\in \mathcal C$, and let $F^n(t)$ be the corresponding renewal-scaled solutions. Suppose $F_t$ is  renewal-scaled solution for initial data $F_0 \in \mathcal P(\mathbb R^+)$ and holding time measure $P\in \mathcal P(\mathbb R^+)$. Then if $u_0^n \rightarrow F_0$ and $p^n \rightarrow P$ weakly, then  $F^n_t \rightarrow F_t$ weakly at all points of continuity of $A^\dagger( t)$.

(2) For any $P, F_0 \in \mathcal P(\mathbb R^+)$, $F_t$ is a probability measure for all $t \ge 0$.

(3) If  $A^\dagger(t)$ is strictly increasing, then for $t>0$,  the map $t \mapsto F_t$ is continuous in $\mathcal P(\mathbb R^+)$ under the weak topology.  \end{theorem}

\begin{remark}
Since $\mathcal C$ is dense in $\mathcal P(\mathbb R^+)$, for any $F_0, P \in \mathcal P(\mathbb R^+)$ we can always find $u^n, p^n \in \mathcal C$ with  $u^n\rightarrow F_0$ and $p^n\rightarrow P$ weakly. Thus, for points of continuity of $A^\dagger$, we can define renewal-scaled solutions through weak limits of strong solutions.  
\end{remark}

\begin{remark}
The possible nonconvergence at jump points  of $A^\dagger$  is essentially due to  the multiple ways that one can define a generalized inverse for increasing functions. In our definition, inverses of cadlag functions remain cadlag, whereas using the traditional definition of a quantile from statistics, for instance, would instead give caglad functions.
\end{remark}

\begin{remark} We note that for statement (3), there are a variety of sufficient $P$ and $F_0$ which produce strictly increasing $A^\dagger(t)$ for $t>0$. One such sufficient condition, for instance, is  if $P(x)$ is strictly increasing on an interval about the origin. The point here is that jumps in the solution only occur when $A(\tau)$ is  constant  in an interval.  As we will see in  Section \ref{ptmass}, point masses arriving at the origin (corresponding to jumps in $A(\tau)$) are continuously redistributed in the $t$ time scale.  
\end{remark}   

\begin{proof}
 \textit{(1)}  Renewal-scaled solutions for  $(F^n_0, P^n)$ satisfy
\begin{equation}\label{ptmeassol}
\bar F_t^n(q) = (1-\bar p^n(q))\int_t^\infty e^{q(\tau^n(t)-\tau^n(s))}ds.
\end{equation}

To show weak convergence, it is enough to show $\bar F_t^n(q)\rightarrow \bar F_t(q) $.  We first wish to show that $A^n\rightarrow A$ weakly.   This follows from the continuity theorem  \cite[XIII.1]{feller2008introduction},  which states that weak convergence
of locally finite measures is equivalent to pointwise convergence of  their
Laplace transforms. From (\ref{fpreinv}),\begin{equation}
\bar A^n(q) = \frac{\bar F^n_0(q)}{1-\bar P^n(q)}\rightarrow \frac{\bar F_0(q)}{1-\bar P(q)} = \bar A(q).
\end{equation}

Weak equivalence of measures implies that cumulative functions $A^n(t)$ converge to $A(t)$ at all continuity points of $A(t)$.  Now, for all points of continuity of $\tau(t) = A^{\dagger}(t)$, $\lim_{n\rightarrow \infty} (A^n)^{\dagger}(t) = A^\dagger(t)$ \cite[Lemma 3.1]{menon2010dynamics}.  As the generalized inverse is a nondecreasing function, its discontinuity set is countable, and thus $\tau^n(t)\rightarrow \tau(t)$ almost everywhere. 

What remains is to show convergence of the integrals

\begin{equation}\label{ptint}
\int_t^\infty e^{-q\tau^n(s)}ds \rightarrow \int_t^\infty e^{-q\tau(s)}ds. 
\end{equation}
Here, we will use the dominated convergence theorem. This requires bounding the integrand $e^{-q\tau^n(s)}$, which we obtain   from several steps:

(i) Observe that there exist $M >0$ and a positive integer $N$ such that $M$ is a point of continuity for $P$, and   $P^n(M)>\frac 12$ for $n>N$. Define truncated holding times which are restricted to $[0,M]$, with cumulative functions 
\be
\tilde P^n(x) = P^n(x)/P^n(M), \quad \tilde P(x) = P(x)/P(M), \quad x \in [0,M].
\ee
Denote  $ \tilde {\mathbf{P}}^n$ and  $\tilde {\mathbf{P}}$ as random variables distributed with respect to $ \tilde P^n$ and $\tilde P(x)$, respectively. Then clearly  $ \tilde {\mathbf{P}}^n \rightarrow \tilde {\mathbf{P}}$ weakly, and since $ P^n(x) \le \tilde
P^n(x)  $, it follows that  $A^n(t) \le \tilde A^n(t).$ 

(ii) Since $\tilde P^n$ is supported on a finite interval,  it has a finite first and second moment. Thus allows us to apply Lorden's estimate \cite{lorden1970excess} to obtain\begin{equation}
\tilde A^n(\tau) \le  \frac \tau{\mathbb E\left[ \tilde {\mathbf{P}}^n\right]}+\frac{\mathbb E\left[( \tilde {\mathbf{P}}^n)^2\right]}{\mathbb E\left[ \tilde {\mathbf{P}}^n\right]^2}+1 . \label{lorden}
\end{equation}
Note that we have added  1 to Lorden's original estimate to account for a possibly distinct initial  holding time.
 
(iii)
Since $ \tilde {\mathbf{P}}^n$ are all supported in  $[0,M]$, it follows that weak convergence to $\tilde {\mathbf{P}}$ implies the convergence of moments. From this and  (\ref{lorden}) it is straightforward to show  that there exist $a,b>0$ such that for all positive  $n$,  
\be
A^n(\tau) \le \tilde A^n(\tau) \le a\tau+b  
\ee
  \begin{equation}
 \Rightarrow \tau^n(t)\ge \frac{t-b}a.
\end{equation}
Thus, integrands in the left hand side of  (\ref{ptint}) are dominated by an exponentially decaying function,  which proves (\ref{ptint}) and thus part \textit{(1)}.

Part \textit{(2)}  follows immediately from the weak convergence shown in part \textit{(1)}, since $F^n$ are all probability measures.

For part \textit{(3)}, note that  if $A(\tau)$ is strictly increasing, then $\tau(t)$ is continuous. From (\ref{measform}),  it is clear that $\bar F_t(q)$ is continuous in $t$, which in turn implies  that $ F_t$ evolves continuously under the weak topology in the space of measures.
\end{proof}
\subsection{An example of  renewal-scaled solutions with  point mass initial conditions}\label{ptmass}

To illustrate how renewal-scaled solutions behave under  jumps in initial data, we conclude with an   example with monodisperse initial conditions.
   Specifically, let the  initial
measure data satisfy $F_0 = \mathbf 1_{[1,\infty)}(x)$ and $p(x) = e^{-x}$.    Then $\mathcal{L}(p(x)) = \frac{1}{q+1},\bar
F_0(q) = e^{-q}$, and   using (\ref{fpreinv}),

\begin{equation}
\bar A(q) = \frac{q+1}{q}e^{-q} = e^{-q}+\frac{e^{-q}}q.
\end{equation}
For $t>0$, Laplace inversion then gives 
\begin{equation}
A(\tau) =\begin{cases}0 & \tau \in [0,1) \\
\tau & \tau \ge 1 \\
\end{cases} , 
\end{equation}
which has a generalized inverse of
\begin{equation}
\tau(t) = A^{\dagger}(t) = \begin{cases}1 & t \in (0,1) \\
t & t>1. \\
\end{cases}
\end{equation}
Substitution of these variables into   (\ref{measform}),  for $t \in (0,1]
$, yields
\begin{align}
\bar F_t(q) = e^{q}\left(\frac {q}{q+1} \int _t^\infty e^{-q\tau(s)}ds\right)
 = e^{q}\left(\frac {q}{q+1}\right)\left((1-t)e^{-q}+ \int _1^\infty e^{-qs}ds\right)
\nonumber
  = \frac{t}{q+1}+(1-t).\nonumber
\end{align}
Taking inverse Laplace transforms then gives us the solution  $F_t(dx) =
(1-t)\delta_0(dx)+te^{-x}dx$ for $t \in (0,1]$.  By a similar calculation, we can show that
$F_t(dx) = e^{-x}dx$ for $t>1$.  

In the $t$ time scale, initial conditions with a support not containing the origin will immediately jump through translation. For our example, the point mass jumps to the origin, and is then continuously reassigned 
  to $[0,\infty)$ according to $p(x)$ at a constant rate of one.  After the reassignment of the
entire point mass, the solution, now a  density, is stationary, since $u(x,\tau) = e^{-x}$
is a solution to (\ref{pde}) with $p(x) = u_0(x) = e^{-x}$ (one can see, in fact, that this is a manifestation of the memoryless property of exponential distributions).

\section{Acknowledgements}The author is indebted to Robert Pego and Govind Menon for multiple conversions regarding the preparation of his thesis~\cite{JK}, a chapter of which    is the basis of this paper. 

\section{Conflicts of Interest}
The author declares no conflict of interest.

\bibliographystyle{siam}
\bibliography{reference2}

\begin{thebibliography}{10}

\bibitem{cox1967renewal}
{\sc D.~R. Cox}, {\em Renewal theory}, vol.~1, Methuen London, 1967.

\bibitem{cox2017theory}
\leavevmode\vrule height 2pt depth -1.6pt width 23pt, {\em The theory of
  stochastic processes}, Routledge, 2017.

\bibitem{feller2008introduction}
{\sc W.~Feller}, {\em An introduction to probability theory and its
  applications}, vol.~2, John Wiley \& Sons, 2008.

\bibitem{kim2004markovian}
{\sc J.~W. Kim and G.~C. Shim}, {\em A {M}arkovian approach to the forward
  recurrence time in the renewal process}, JKSS (Journal of the Korean
  Statistical Society), 33 (2004), pp.~299--302.

\bibitem{JK}
{\sc J.~Klobusicky}, {\em Kinetic limits of piecewise deterministic {M}arkov
  processes and grain boundary coarsening}, PhD thesis, Brown University,
  Providence, RI, 2014.

\bibitem{liao2013applied}
{\sc M.~Liao}, {\em Applied stochastic processes}, CRC Press, 2013.

\bibitem{lorden1970excess}
{\sc G.~Lorden}, {\em On excess over the boundary}, The Annals of Mathematical
  Statistics,  (1970), pp.~520--527.

\bibitem{menon2010dynamics}
{\sc G.~Menon, B.~Niethammer, and R.~Pego}, {\em Dynamics and self-similarity
  in min-driven clustering}, Transactions of the American Mathematical Society,
  362 (2010), pp.~6591--6618.

\bibitem{resnick2013adventures}
{\sc S.~I. Resnick}, {\em Adventures in stochastic processes}, Springer Science
  \& Business Media, 2013.

\bibitem{stewart2009probability}
{\sc W.~J. Stewart}, {\em Probability, Markov chains, queues, and simulation:
  the mathematical basis of performance modeling}, Princeton University Press,
  2009.

\end{thebibliography}

\end{document}